\input amstex
\documentstyle{amsppt}
\input epsf
\input psfrag

\define\RR{\Bbb R}
\define\CC{\Bbb C}
\define\SS{\Bbb S}
\define\FF{\Bbb F}
\define\ZZ{\Bbb Z}
\define\DD{\Bbb D}
\def\symboli{\text{\bf i}\,}
\def\symbolj{\text{\bf j}\,}

\topmatter
\title Complex odd-dimensional endomorphism and topological degree
\endtitle
\rightheadtext {}
\author Jon A. Sjogren
\endauthor
\affil Faculty of Graduate Studies, Towson University
\endaffil
\address Towson, Maryland
\endaddress
\date 1 September 2015
\enddate
\endtopmatter
\document

\abovedisplayskip=13pt
\belowdisplayskip=13pt

\head Abstract\endhead

Any linear mapping on the real vector space $\RR^n$, where the dimension $n = 2m+1$ is odd, should have a real eigenvalue-eigenvector pair. It is also true that a linear $T: \CC^n \to \CC^n$ has a complex eigenvector  $v \in \CC^n$, $T v = \mu v$, where the eigenvalue $\mu \in \CC$, regardless of the divisibility of $n \geq 1$. But here again, when $n = 2m+1$, the result should be more natural and easier to prove. H. Derksen uses this result as the base case for an induction leading to ``$T$ has a complex eigenvector'' for any $n$ at all (the Fundamental Theorem of Algebra). We build up enough of the homotopical theory of $k$-plane bundles on a single base space, $B = \RR P^2$, to obtain ``Derksen's Lemma'', meaning the case where $n = 2m+1$. Two technical results are required which are proven completely within 3-dimensional Euclidean topology (including the {\it degree} of a mapping $h: \SS^2 \to \SS^2$). Let $\gamma$ be the canonical (``Hopf'') line bundle on $B$, and $\epsilon$ the trivial line bundle. Firstly, $\gamma \oplus \gamma \oplus \gamma \oplus \gamma \simeq \epsilon \oplus \epsilon \oplus \epsilon \oplus \epsilon$. In other words, $4\gamma$ (Whitney sum) is a trivial bundle of rank $4$. Secondly, $\gamma \oplus \gamma \oplus \epsilon$ is {\it not} trivial, it is not $\RR P^2$-isomorphic to $\epsilon \oplus \epsilon \oplus \epsilon$.

\abovedisplayskip=13pt
\belowdisplayskip=13pt

\head Synopsis\endhead

H. Derksen has produced an interesting proof of the Fundamental Theorem of Algebra that relies on properties of certain operators on linear spaces of symmetric matrices over the complex numbers. The ``analytic'' part of his proof seems to lie in the well-known fact that an odd-order polynomial over the {\it reals} $\RR$ does have a root. Using the companion matrix construction, Derksen poses the problem as finding an eigenvector for any finite endomorphism $T: \CC^n  \to \CC^n$. A proof by induction proceeds according to the size of the integer $\nu$ where $n = 2^{\nu} (2m+1)$, just as in the proof by P.S. Laplace (1795), and later proofs due to Gauss, Gordan and E. Artin! A reference for the history of the FTA is \cite{Fine \& Rosenberger}.

Remarkably, the base case of Derksen's induction, $\nu = 0$, uses a different matrix representation than does the inductive step. It could be convenient to have a separate proof of this 0-level result, one that is consistent with the ``linear algebra'' framework. Simply put, we need to exhibit an eigenvector for $T$ when $n$ is odd.

Toward this goal, one works with parametrized linear algebra over the real projective plane $\RR P^2$ (and sometimes over the projective line). The degree of the original polynomial is revealed in the rank of a pair of Whitney-summed hyperplane bundles: the trivial $2n$-bundle, and a $2n$-fold sum of copies of the canonical (``Hopf'') line bundle. Given some operator $T_0$ on $\CC^n$ with no eigenvalue-eigenvector, results from the 1980s and 90s, on the topology of matrix spaces, show that these must be $B$-isomorphic bundles, with base space $B = \RR P^2$. This improbable situation can be repudiated by the application of the (total) Stiefel-Whitney class. However, our intention in this article is to operate within a limited range of mathematics. This range certainly includes the definitions and basic facts of the homotopy concept, including its application to line bundle constructions
over a compact Hausdorff space.

We keep the reader informed of alternate proofs involving various methodologies. But in principle we avoid homology groups, cohomology rings and their operations, characteristic classes, K-theories with computation by spectral sequence, classifying spaces (and their cellular architecture!), as well as techniques from differential topology, such as critical points, transversality, and ``pre-image'' manifolds. Simplicial approximation, measure and density results are held to a minimum.

What does enter in are propositions pertaining to Brouwer degree for dimension two only. We refer to the ``classical'' Borsuk-Ulam theorem which is an introductory result of geometric topology. How to express the topological degree (for maps on $\SS^2$) is handled according to preference, either through linear approximation, differential forms such as ``curvatura integra'', or even by counting pre-image points.

We hope that the one or two new technical  observations made here will  be applicable outside the limited domain where we have used them. Whether or not some truly new results might arise, it is worthwhile to seek the underlying reason as to why an ``odd transformation'' of a complex space has an ``axis'', as does a spinning spheroid in mundane 3-space. 

\head Introduction\endhead

Solving a polynomial equation (in one variable) has long been known to the tied to the ``eigenvector problem'' for a linear mapping. Given  $T: \FF^n \to \FF^n$ linear, where $\FF$ is a field, the characteristic (or secular) polynomial is  defined as $P(X) = \det \left[IX -T\right]$ where a  matrix realization of $T$ has been chosen, $X$ being  the indeterminate generator of $\FF[X]$. Given $P(X) = X^n + \alpha_{n-1} X^{n-1} + \cdots + \alpha_0$, the ``companion matrix'' ${\Cal C}_p$ yields an endomorphism with  respect to a certain basis, which returns $P$ as its characteristic polynomial (irrespective of the basis).

This correspondence is used to solve problems in operator theory using polynomial algebra, and vice-versa.
We discuss two known cases of this correspondence phenomenon, one of them rather more prototypical. Taking $\FF = \RR$, the standard ``real numbers'', it has been known since antiquity  that  rigid motion of a ball in three-space possesses an axis of rotation, while in the plane there generally is none.

The correspondence poses the question of  $P(X) = X^n + \beta_{n-1} X^{n-1} + \cdots  = 0$ having a root (``eigenvalue of $T$'') where $n$ is {\it odd} and $\beta_j \in \RR$.
This fact follows from the Intermediate Value Theorem, or one may say, the topological closure of the reals $\RR$.

Although a linear analysis can  be done  keeping use of ``determinant'' to a minimum, see \cite{Axler}, one might expect to prove the above ``Axis Theorem'' entirely without polynomial algebra. We point out two ways to do this, both coming  from topology.

We now step ahead to explain that we are not only engaged in re-proving well-known facts. To take $\FF = \CC$, the complexes, we may ask (where determinants as a technique are off-limits) whether such an operator $T : \CC^n  \to \CC^n$, $n$ odd, need have any ``axis'', that is, a one-dimensional subspace, or eigenvector. If this is proven so, the  companion correspondence shows that the Fundamental Theorem of Algebra for a complex polynomial $P(z)$ of degree $n = 2m+1$ is true, i.e., $P$ has a root $z_0 \in \CC$.

It is well-understood that this ``$n$ odd'' case of FTA is easier than the general case. One class of FTA proofs operates by reducing  the given $n$ by factors of two in the $\RR$ case. For such a proof (de Foncenex, Laplace 1795, Gauss 1816, Gordan et al) only one reduction step, instead of many, is required. Leaping  to another
category, FTA (odd over $\CC$) is proved earlier on in the book \cite{Guillemin \& Pollack} than is the general FTA. The former uses, from the point of view of differential topology, the  relatively primitive concept of ``mod 2 intersection (or degree)'', whereas the latter requires full integer intersection or degree invariants.

We have already alluded to topological methods. At this point we admit that the  ``full'' FTA (or ``every finite $\CC$-linear mapping has an eigenvector'') now has a short topological proof due to \cite{Sohsten de Medeiros}. The Medeiros proof uses the Lefschetz Fixed-Point Theorem applied to a self-mapping  of complex projective space induced by the $\CC$-linear mapping. For another ``real'' proof, based on Brouwer's Invariance of Domain principle, see \cite{Sjogren}.

We conclude this Introduction with two remarks. Firstly, we give as motivation the work of others, mainly \cite{Derksen}. Secondly, we supply a list of mathematical concepts and techniques that will be avoided, so that the ``new proof'' is clear by its prerequisites.

H. Derksen, see also \cite{Conrad}, sets up an induction over $n$ leading to an eigenvector for the  given endomorphism of $\CC^n$. In fact, stronger results are proved, concerning commuting sets of morphisms that turn out to possess an eigenvector in common. The ``base case'' of the induction amounts to our  Axis Theorem for $\CC^{2m+1}$. The Derksen proof  makes the matrix $A$ for $T$ act on a linear space of Hermitian matrices in two different ways.

The two actions happen to commute, which  allows use of a general result due to Derksen on commuting operators. One point is that we have not proven FTA as yet, since  only the base at induction (our $2m+1$ case) has been done. The induction step involves different linear spaces than before. In \cite{Conrad} it is remarked that the base case of an induction is  not necessarily a trivial proposition.
It might therefore seem desirable to  provide a separate proof of this ``base'', whence it can be referred to as needed. The new ``base'' proof should be reasonably consistent in  methodology to Derksen's concept.

We arrive at the disclaimer portion of  the Introduction. First we assert that there are two approaches to both the ``odd'' ($n = 2m+1$) eigenvector case for $\RR^n$ and for $\CC^n$. These we named the ``vector field method'' and the ``line bundle method''. Each ``method'' has in  turn two separate approaches, which can be called the ``mechanical'' and ``handmade'' technologies.

As a first example, we might want to  show that a certain class of $k$-manifold $\{M\}$ has only tangent vector fields that vanish at some point $y \in M$. The ``mechanical'' technology would  employ say cohomology theory, which brings algebra to bear. A ``handmade'' proof  might start with a non-vanishing field on $M$ and arrive at a geometric incongruity, such  as a continuous function on an interval that misses an intermediate value.

Our preferred proof is the handmade ``line bundle'' version. We work in a linear category, but ``parametrized'' over a manifold. The base space required is either  $\RR P^1 \simeq {\SS}^1$ or $\RR P^2$, where the construction of vector bundles is well-understood. The ``mechanical'' approach provides an alternative, employing the naturality of  characteristic classes. On the other hand, the  recommended proof avoids all homology groups,  cohomology (graded) rings, Steenrod algebra, K-rings and associated spectral sequences, characteristic classes in ordinary and extraordinary cohomology, and techniques typical of differential topology including regularity, transversality, and pre-image manifolds. Measure, probability and simplicial approximation play a minimal role. The tools that  are exploited have a firm geometric foundation such as Brouwer's Fixed-Point Theorem. The Borsuk-Ulam Theorem is used only in its 2-dimensional formulation which has a simple general topology proof. Elementary ``degree theory'' (mappings of a k-sphere $\SS^n$) comes in, where the reader has a choice of a basic simplicial  argument (that homotopy implies equal degree),  or a Kronecker integral type construction due to  Mawhin-Heinz. In fact we admit that  the theorem of Sard-Brown for $C^1$ functions on $\SS^n$ is rather straightforward, so a ``pre-image'' formulation of topological degree is not far-fetched.  The underpinnings  of our approach lie in the ``elementary'' homotopy theory of vector bundles on a compact base.

The present (linear) theory does not rely upon results from general  (principal) fiber bundles, but the subject of line bundles up to homotopy has not been widely exposited. We do review its major suppositions and consequents below, see also \cite{Husem\"oller}.

\head The Method of Vector Fields\endhead

Let $n = 2m+1 > 0$ be odd. Then a linear transformation $T_r: \RR^n \to \RR^n$ or $T_c: \CC^n \to \CC^n$ gives rise to a vector field on the sphere $\SS^{n-1}$ or to two vector fields on $\SS^{2n-1}$. In the  real case, if $T_r$ has no real eigenvalue/eigenvector,  the vector field has no vanishing point (it is  ``independent'' throughout $\SS^{2n-1}$). Similarly if $T_c$ has no complex eigenvalue/eigenvector, the two fields also constitute an independent set of fields on $\SS^{2n-1}$. The fact that such collections of  independent vector fields can be shown topologically not to exist means that $T_r$ and $T_c$ (in dimension $n$) have their eigenvectors. Still, not all continuous  or smooth vector fields arise from this kind of ``matrix'' construction: we used a tool more powerful than really needed. It may be asserted that the ``line bundle'' approach later to  follow, is better adapted to our $2m+1$ eigenvector problem.

We move on to the ``real, odd'' situation. Let $T$ be linear: $\RR^{2m+1} \to \RR^{2m+1}$, then if $T$ is singular, it has $\lambda = 0$ as an eigenvalue and an associated characteristic vector or eigenvector $v_0$ so  that $T v_0 = \vec{0} \in \RR^n$. If $T$ is non-singular,  it can be restricted as a continuous mapping  to $\SS^{2m} \subset \RR^n$, leading to $A(x) = \frac{T(x)}{||T(x)||}$, scaled by the Euclidean norm.
If $T$ has no eigenvector, then $A: \SS^{2m} \to \SS^{2m}$ never maps any $y \in \SS^{2m}$ to $y$ (itself) nor to  its opposite (antipode) $-y$. But then, a homotopy can readily be constructed between $\text{id}_{\,\SS^{2m}}$ and the antipodal mapping  $y \longmapsto -y$, see \cite{Dugundji}. Since the  dimension $2m$ is even, these mappings have Brouwer degree $= + 1$, respectively $-1$, so cannot be freely homotopic. Hence an eigenvector $v$ of $T$, with its real eigenvalue $\lambda$ must  exist. Essentially the same ``mechanical'' proof can be carried out by looking at the  mapping that $A$ induces on the integral homology group $H_{2m}$ ($\SS^{n-1}$).

For a handmade proof one should appeal to the famous ``cue-ball'' theorem as proved in \cite{Milnor CUE} using ideas from foliation theory of D. Asimov. No ``invariants'' are used, only the  calculus of volumes and the Contraction Mapping Principle. The vector field $A(s)$ defined everywhere on $\SS^{2m}$, taking values in $\RR^n$, can  be flattened smoothly to a tangent vector field on $\SS^{2m}$, one that never vanishes. But it is well-known that such a tangent vector field does not exist, as shown for  instance by the Milnor-Asimov proof.

Thus we may be satisfied with the ``real, odd'' situation. We found an ``axis'' for any linear transformation $T$ on $\RR^{2m+2}$. Furthermore, using the correspondence between  real matrices brought in by the ``companion'' construction, we arrive at the profound result that
$$X^{2m+1} + \alpha_{2m} X^{2m}  + \cdots + \alpha_0 = 0$$
always has a real root, without our needing to  form any kind of analytic estimate using the coefficients themselves; compare \cite{Hille, p. 208}.

\head Spanned Matrices of Constant Rank\endhead

The problem we just now addressed, existence of a real eigenvector, is part of a theory dealing with  vector spaces that consist of matrices whose rank is fixed, with the square matrix of  ``all zeros'' thrown in. We could have asked, ``are there nonsingular matrices $M_{2m+1}$ and $N_{2m+1}$, neither a scalar multiple of the other, that {\it span} nonsingular matrices, together with the  zero matrix?'' We found that the answer was ``no''. Supposing that $n$ is not necessarily odd but rather $n = 2^b (2m+1)$, where $b = c + 4d$, $0 \leq c < 4$, we have the Theorem of Radon-Eckmann-Adams: \proclaim{Invertible Spanning} The size of a maximal independent set of invertible $n\times n$ matrices $\{A^j\}$ over $\RR$, whose span with coefficients $\{\lambda_j\}$, not all zero, is  $\rho(n) = 2^c + 8d$, the value of the Radon-Hurwitz  function at $n$.\endproclaim See \cite{Adams, Lax, Phillips} and related subsequent literature.

We confirmed this result above when $b = 0$. We are going to go over this instance once more, consistent with the geometric innovation proposed here. After that, the  article concerns the case where $b = 1$. We are seeking  geometric methods that generalize, but apparently are not related to ``$\psi$-operations in K-theory'' or spectral sequence calculations arising from elaborate cellular decompositions. Let us examine how the ``vector field approach'' also lends itself to the case $b = 1$.

\head The Vector Field Approach to Complex Eigenvector (in Odd Dimension)\endhead

We briefly cover the general case of $T: \RR^n \to \RR^n$ where $n$ can be any integer $\geq 1$. Suppose that $T_1, \dotsc, T_p$ are nonsingular. If $T_k^{-1}T_j v = \lambda v$ for a non-zero vector $v$, it would mean that some non-trivial linear combination of the $\{T_i\}$ turns out {\it not} to have full rank. Otherwise, $T_k^{-1}T_j$ restricted to $\SS^{n-1}$ gives rise to a vector field without ``poles''; it can everywhere be flattened to a unit vector field on $\SS^n$, see \cite{Milnor TFDV}.

On the other hand, if for some $u \in \SS^n$, $\displaystyle{\sum_{i=1}^{k-1}}\, \beta_i A_k^{-1} A_i u = 0$ as a tangent vector, it would mean that some non-trivial linear combination of matrices $\sum \gamma_i T_i - T_k$ has a null vector, violating the stipulation that such a non-trivial combination has full rank, see \cite{Causin}. Thus, restrictions on the size of a set of independent vector fields tell you also how many nonsingular ``spanning'' matrices  can exist in the given dimension. Similar techniques have been used to study spanning sets of symmetric and Hermitian matrices, and of real matrices whose span consists of matrices not invertible, but of some  fixed rank. See \cite{Petrovi\'c}, \cite{Lam \& Yiu} and \cite{Meshulam}.

We step back to the humble case of our present interest. Given $T$ an $n\times n$ matrix over  $\CC$ where $n$ is odd, we cast everything in terms of real vector spaces and their transformations. Rewriting $T = \left[c_{ij}\right], c_{ij} \in \CC$ by $A = \left[\matrix
a_{ij} & -b_{ij} \\
b_{ij} & a_{ij}
\endmatrix\right]\,$, where $c_{ij} = a_{ij} + {\roman i} b_{ij}$,  $a_{ij}, b_{ij} \in \RR$, we obtain a real matrix acting on $\RR^{2n}$.
If $A$ has a real eigenvector $v \in \RR^{2n}$, then $T$ has an eigenvector $\left(v_1 + {\roman i}v_2, v_3 + {\roman i} v_4, \dotsc, v_{2n-1} + {\roman i} v_{2n}\right)$ with {\it real} eigenvalue $\lambda$.

If $T$ has no complex eigenvector, then $A = A_T$ and also ${\roman i} A = A_{\,{\roman i}T}$ have no real eigenvectors, and furthermore span a space of nonsingular matrices. Indeed, if $\left(\beta + {\roman i}\gamma\right) A$ is  singular, then $A$ has a complex null vector, which already yields a contradiction. Moreover we cannot admit $\alpha  I_{2n \times 2n} + \beta A + {\roman i}\gamma A$ as singular either, since this would imply that there exists a complex vector $w \in \CC^n$ such that $T w = \rho w$, where 
$$\rho = \dfrac{-\alpha}{\beta + {\roman i} \gamma}\ \,\,.$$  We have used notation to facilitate the exposition, such as $T \sim A = \left[a + {\roman i} b\right] = \left[\matrix
a & -b \\
b & a
\endmatrix\right]$, ${\roman i}A = \left[\matrix
-b & -a \\
a & -b
\endmatrix\right]$ and so forth.

We established that a complex matrix $T$ of odd square dimension $n = 2m+1$, with no eigenvector, leads to  matrices $\{I, D, E\}$, real of size $2n \times 2n$ that span only nonsingular matrices and the 0 matrix. But this outcome is seen to be impossible by results of Eckmann, G. Whitehead,  J.H.C. Whitehead and  Steenrod, Adams and others.

We intend to sketch out a  classical proof, and then offer a new proof based on ``parametrized linear algebra'' where the parametrizing space is just the real projective plane $\RR P^2$. We concluded above that for $n = 2m+1$, a complex $n \times n$ matrix has an  eigenvalue-eigenvector pair, or else there is a  triple of nonsingular $2n \times 2n$ real matrices spanning only invertible matrices (and the 0). We also saw that this situation leads to  the construction of two independent vector fields on $\SS^{2n-1} = \SS^{4m+1}$. The mentioned classical results in homotopy theory rule this out.

At about the same time, B. Eckmann in Switzerland and G. Whitehead at Chicago gave proofs that the $4m+1$-sphere has one independent (non-vanishing) tangent vector field but no two such fields. We call these proofs our ``handmade'' derivations because they deal with concrete constructions, mostly on some space of linea subspaces (Grassmann) or a space of frames  (Stiefel), in a limited range of dimension (case $b > 1$).

There remains to summarize the ``mechanized'' approach to fields on $\SS^{4m+1}$. Since the dimension is odd, placing the sphere into $\RR^{2n}$ in the natural way affords one vector field, namely
$$v(S) = v\left(x_1, x_2, \dotsc, x_{2n-1}, x_{2n}\right) = \left(-x_2, x_1, \dotsc, -x_{2n}, x_{2n-1}\right)$$
for $s \in \SS^{2n-1} \subset \RR^{2n}$.

Keeping in mind the possibility of Gram-Schmidt orthonormalization of two independent fields, it is enough to show that an ``orthonormal pair'' of vector fields on $\SS^{2n-1}$ does not exist, see \cite{Whitehead-Steenrod}. In general, given the space of $k$-frames $V_{n,k}$ in $\RR^q$, $q = n$ or $q = 2n$, there is a fibration $\rho = V_{q,k} \to \SS^{q-1}$.

Thus a {\it section} $f: \SS^{q-1} \to V_{q,k}$ yields  a $k$-frame (orthonormal) for each point $s \in \SS^{q-1}$. The assumption that there are two independent vector fields causes us to consider the sections of $V_{q,3}$, where now $q = 2(2m+1)$. The section $f$ induces a surjection of  cohomology classes
$$f^*: H^{q-1} (V_{q,3}) \to H^{q-1} (\SS^{q-1})\,\,.$$
\indent
Our smallest case of interest is represented by $q = 6$. In general, if $2k-1 \leq q$, the $(q-1)$-skeleton of $V_{q,k}$ looks like
$$\RR P^{q-1} /\RR P^{q-k-1} = \wp_{q,k}\,.$$
\indent
This latter space is explicitly understood from its cell structure as a CW-complex. We have a short exact sequence of mappings
$$(a) \to \RR P^{q-k-q} \overset{\symboli}\to{\longrightarrow} \RR P^{q-1} \overset{\pi}\to{\longrightarrow} \wp_{q,k} \to (b)\,, \tag{*}$$
where $a,b$ are base points. It is well-known, \cite{Hatcher AT}, that
$$H^* \left(\RR P^{q-1}, \ZZ_2\right) \simeq \ZZ_{2} [\alpha]/(\alpha^q)\,\,,$$
where $\ZZ_2$ is the commutative ring with two elements, zero and unity, and $(\alpha^q)$ is the principal ideal of the polynomial ring indicated, so also
$$H^* \left(\RR P^{q-k-1}, \ZZ_2\right) \simeq \ZZ_2 [\beta]/\left(\beta^{q-k}\right)\,,$$
where ${\symboli}^* \alpha = \beta$.

From the short exact sequence (*) and its long exact cohomology sequence, we deduce
$$H^{q-1} \left(V_{q,k}, \ZZ_2\right) \simeq H^{q-1} \left(\wp_{q,k}, \ZZ_2\right) \simeq \ZZ_2\,,$$
since both spaces look the same up to the $q$-skeleton.

As a $\ZZ_2$-module, we see that $H^{q-1}(V_{q,k})$ has a generator $w_{q-1}$, also $H^{q-1}(\SS^{q-1})$ has a  generator $\sigma_{q-1}$. As mappings we have $\rho \circ f = \text{id}_{\SS^{q-1}}$, so passing to $(q-1)$-cohomology, we must have $\rho^* \sigma = w$ and $f^* w = \sigma$. Specializing to $k = 3$, there is a Steenrod formula for $H^* \left(\RR P^{q-1}, \ZZ_2\right)$ as a ring (with cup product):
$$S q^2 \alpha^q = \binom{q}{2} \alpha^{q+2}\ .
$$
\indent
This leads to a formula for $v_{q-3}$, the additive generator of $H^{q-3} \left(\wp_{q,3}, \ZZ_2\right)$, namely
$$S q^2 v_{q-3} = \binom{q-3}{2} v_{q-1}\,,$$
where  $v_{q-1}$ generates $H^{q-1} \left(\wp_{q,3}, \ZZ_2\right)$. Since $q = 2n$, it is clear that
$$\binom{q-3}{2} = \binom{2(2m+1)-3}{2} = (4m-1)(2m-1)\,,$$
the product of two odds.

Hence over $\ZZ_2$ we have $S q^2 \,v_{2n-3} = v_{2n-1}$. By the naturality of squares, $\sigma = f^* w = f^* S q^2 v_{2n-3} = S q^2 f^* v_{2n-3}$. But $f^* v_{2n-3} \in H^{2n-3} \left( \SS^{2n-1}, \ZZ_2\right)$. This latter $\ZZ_2$-module is the zero module if  $2n> 3$, so we have a contradiction inasmuch as $\sigma$, the generator of $H^{2n-1}\left(\SS^{2n-1}\right)$, equals $S q^2 (0) = 0$. Hence the desired section $f: \SS^{2n-1} \to V_{2n,3}$ cannot exist, nor can a continuous 2-field on $\SS^{4m+1}$.

The preceding ``vector field approach'' (mechanical version) is based on a sophisticated theory, firstly bearing on the geometric structure of ``flag manifolds'', and secondly related to the formal algebra of  cohomology operations. In addition, the way to a proof (but also when $b = 2,\, b=3, \dots$) calls upon a fair amount of technique. Working instead with line bundles requires a foundation only in the elementary homotopy of  vector bundles, which provides effective insight without getting overly algebraic. The spaces we work with have dimension similar to that of the spanning set, in our case of interest, $k = 3$, {\it not} the size of the original real dimension, $q = 2n$. The relevant bundles can to an extent be visualized. Thus their triviality or non-triviality can sometimes be apprehended directly .

\head Line Bundle Method for $b = 0$ and $b=1$\endhead

Elementary definitions on the category of real vector  bundles are taken for granted. The basic propositions  are reviewed; the reader should have access to \cite{Hatcher VB}, \cite{Milnor \& Stasheff}, \cite{Husem\"oller}, and
\cite{Dugger}. Two thrusts of the elementary homotopy theory of vector bundles (not referring to general fiber bundles) are
\roster
\item"(i)"
establishment of real linear algebra constructions smoothly over a parameter space (the base) and
\item"(ii)"
the fact that the {\it dimension} of the base limits
the discrepancy between the concepts ``stably trivial vector bundle'' and ``trivial vector bundle''.\endroster
If in some case a more general fiber construction comes up, it is likely to be a homotopy sphere such as $F = \RR^p - \{0\}$.
\bigskip
\noindent
{\bf Proposition 1.}  Suppose $B$ is a CW-complex of dimension $m$, and $F$ is the fiber of a bundle $(E, \pi, B)$ where $E$ is the total space and $\pi$ is the projection. Let $A \subset B$ be a subcomplex, then given a section $s: A \to E$,  it may be {\it extended} to a ``global'' section $s^*: B \to E$ in case $F$ is $(m-1)$-connected.

In other words, if for $1 \leq k < m$, any mapping of $\SS^k$ into $F$ is contractible to a point within $F$, the section $s_{|A}$ can be ``prolonged'' to a valid section on $B$. This proposition boils down to the fact that
the ``contractible'' assumption implies that an attaching map $\delta K \to F$ in a cellular construction extends continuously to the interior of the cell $K$.

We assume that our base space is a finite CW-complex hence paracompact (any open cover gives rise to a subordinate partition of unity). In fact, the only base we need in our application is a compact manifold, hence paracompact \cite{Hatcher VB}.

The following is found in all the textbooks.

\bigskip
\noindent
{\bf Proposition 2.} Let $u: \xi \to \xi'$ be a vector bundle morphism over a base $B$. Then $u$ is an isomorphism if and only if $u: \pi^{-1}(b) \to (\pi')^{-1}(b)$ is a vector space isomorphism for each $b \in B$.

Now let $I = [0,1]$, the closed unit interval.

\bigskip
\noindent
{\bf Proposition 3.} Given $\xi = \left(E, \pi, B \times I\right)$ with paracompact base $B$, there is a bundle isomorphism
$\xi | B \times 0 \simeq \xi | B \times 1$, consistent with the retraction $r: B \times  I \to B \times 1$.

This follows from Proposition 2 and the fact that  there is a bundle morphism $(u, r): E \to E$ that is an isomorphism on each fiber.

\bigskip
\noindent
{\bf Proposition 4.} Let $u: \xi^n \to \eta^m$ be a bundle morphism over $B$, whose rank on every fiber $u_b: \pi_{\xi}(b) \to \pi_{\eta}^{-1}(b)$ is an unvarying nonnegative integer $k$. Then, the vector spaces $\text{ker}\,u (b)$, $\text{im}\,u (b)$ and $\text{coker}\,u(b)$ form locally trivial vector bundles with the evident projection map.

This significant result arises from examination of
$$u_b: F^n = V_1 \oplus \text{ker}\,u_b \to F^m = \text{im}\,u \oplus W_2\,.$$
It may be shown that an induced linear mapping
$$w_b: V_1 \oplus \text{ker}\,u_b \oplus W_2 \to \text{im}\,u_b \oplus W_2 \oplus V_2$$
is an isomorphism over an open set $U \subset B$ that  contains $b$. Inverting this map leads to a  (continuous) section $u \mapsto v_b = w_b^{-1}$, which leads to triviality on $U$ for the three proposed bundles of interest.

\bigskip
\noindent
{\bf Proposition 5.} If a $B$-morphism $u:\xi^n \to \eta^m$  is a monomorphism on each fiber of $\xi$ (injective), then $\text{im}\,u$ and $\text{coker}\,u$ are vector bundles. If $u$ is surjective, an epimorphism on each fiber of $\xi$, then $\text{ker}\,u$ is a vector bundle.

\bigskip
\noindent
{\it Proof.} Such $B$-morphisms have constant rank, so we apply Proposition 4.

The following important result is proven using a Riemannian metric on the bundle, or an injective mapping $E \to F$ called the ``Gauss map''. Since any real bundle over a paracompact $B$ has such a metric, and a Gauss map, this Proposition is valid for our cases of interest.

\bigskip
\noindent
{\bf Proposition 6.} Given $0 \to \xi \overset{u}\to{\longrightarrow} \eta \overset{v}\to{\longrightarrow} \zeta \to 0$, a short exact sequence of vector bundles oven a CW-complex $B$, with a Riemannian metric on $\eta$. Then there exists a ``splitting $B$-morphism'' $w: \xi \oplus \eta$, satisfying $w {\symboli} = u$, $\upsilon w = {\symbolj}$, where ${\symboli}: \xi \to \xi \oplus \zeta$ is the natural injection and ${\symbolj}: \xi \oplus \zeta \to \zeta$ is the projection.  

\bigskip
\noindent
{\bf Remarks.} The given metric allows one to build a continuous mapping $g: B \times F^m \to B \times F^n$ which induces an epimorphism $g: \eta \to \text{im}\,u$.

\bigskip
The mapping $w$ on $\xi$ is just the isomorphism $\xi \to \text{im}\,u$, and on $\zeta$, $w$ is defined as mapping to the orthogonal complement of $\text{im}\,u$ in $E(\eta)$ according to the metric. See \cite{Husem\"oller} for fuller details.

We have just covered the main results showing how to do most of linear operator theory over a ``base''. Now suppose that the base $B$ is $d$-dimensional CW-complex or manifold. The proofs of the following use Propositions 1 through 3. See \cite{Dugger, 11.10} for a less ``elementary'' derivation, using classifying spaces. Now let $\epsilon$ be the ``trivial'' or product line bundle over $B$.

\bigskip
\noindent
{\bf Proposition 7.}  If $k \geq d+1$, any $k$-bundle $\xi$ over $B$ is isomorphic to $\eta \oplus \epsilon$ for some $k-1$-bundle $\eta$.

\bigskip
Given this, there follows by induction:

\bigskip
\noindent
{\bf Proposition 8.} If $\xi_1$ and $\xi_2$ are two $k$-bundles with $k \geq d+1$, and $\xi_1 \oplus \epsilon^{\ell} \simeq \xi_2 \oplus \epsilon^{\ell}$, for some $\ell \geq 0$, then $\xi_1 \simeq \xi_2$.

\bigskip
Thus, a real vector bundle of dimension $k$ over a paracompact $B$ of dimension $d$, satisfying $k \geq d+1$, can be written in the form $\eta \oplus \epsilon ^{k-d-1}$ where $\eta$ is a $d+1$-vector bundle which is uniquely determined up to isomorphism. In particular if $\xi^{d+1}$ is stably trivial, i.e. $\xi^{d+1} \oplus \epsilon^{\ell} \simeq \epsilon^{d+1+\ell}$ then $\xi^{d+1} \simeq \epsilon^{d+1}$.

\head Matrix Spaces and  Line Bundles\endhead

Let ${\Cal M}(q, k)$ be the set of real $q \times q$ matrices of rank $k$, and ${\Cal V}(q, k)$ be some $r$-dimensional linear subspace of ${\Cal M}(q, k) \cup \{0\}$. Thus, ${\Cal V}$ is spanned by matrices $A_1, \dotsc, A_r$ for which there is no non-trivial linear combination other than results in another $n\times n$ real matrix of rank $\geq k$. Consider the compact manifold $\RR  P^{r-1}$. This space affords two well-known and distinct line bundles, the canonical (or ``Hopf'') bundle $\gamma$,     and the trivial line bundle $\epsilon$ whose total space is $E_{\epsilon} = \RR P^{r-1} \times \RR$ with the natural projection onto the  base. Keeping in mind the defining projection
$$\pi: \SS^{r-1} \to \RR P^{r-1}\, ,$$
which identifies $s \in \SS^{r-1}$ with $-s$ consistent with the usual embedding $\SS^{r-1} \subset \RR^r$, we form a linear combination of the given real matrices,  based on $s = \left(x_1, \dotsc, x_r\right)$, namely
$$A(s) = x_1 A_1 + \cdots + x_rA_r\,\,.$$
A typical {\it point} in the total space of $q \cdot \gamma$ (sometimes written $\gamma^q = \gamma \oplus \cdots \oplus \gamma$, a $q$-fold Whitney sum) can be expressed as
$$\left(\pm s, \lambda_1 s, \dotsc, \lambda_q s\right)\,,$$
see \cite{Milnor \& Stasheff}, where $\{\lambda_j\}$ are real parameters. We may define
$$A \left(\pm s, \lambda_1 s, \dotsc , \lambda_q s\right) = \left(\pm s, \displaystyle{\sum_{i=1}^r}\,x_i A_i (\lambda_1, \dotsc, \lambda_q)^T\right),$$
which gives a point in the total space of $\epsilon^q$ with fiber $\RR^q$ over the base $\RR P^{r-1}$. By the  ``linear algebra'' propositions above, we observe that  $A$ is then a constant-rank $k$ bundle morphism $\gamma^q \to \epsilon^q$, hence $\text{Im}\,A$ is a vector bundle of dimension $k$. Furthermore, splittings exist of the form
$$\align
\text{Im}\,A \oplus \text{ker}\,A & \simeq q \gamma_{r-1}\,\, , \tag{\dag} \\
\text{Im}\,A \oplus \eta^{q-k} & \simeq \epsilon^q\ ,
\endalign
$$
where $\eta$ is the $(q-k)$-plane bundle complementary to $\text{Im}\,A$, a bundle that uses the Riemannian  metric in its construction. Useful references for the above, and generalizations are \cite{Petrovi\'c, thesis}, \cite{Lam \&  Yiu}, and \cite{Meshulam}.

Our cases of interest concern square matrices of full rank $q$, which can be an odd or twice an  odd integer. Since $\text{ker}\,A$ and $\eta$ are then zero-dimensional bundles, the isomorphisms from  $(\dag)$ collapse and we are left with
$$q \,\gamma_{r-1} \simeq \text{Im}\,A \simeq \epsilon^q\,\,,$$
so the $q$-fold canonical line bundle is actually the trivial $q$-plane bundle.

First we cover the situation of two real square matrices $A_1, A_2$ of size $q = n = 2m+1$,
both invertible, neither a real multiple of the other. We arrive at an isomorphism of vector bundles
$$\underbrace{\gamma_1 \oplus \cdots \oplus \gamma_1}_n \simeq \underbrace{\epsilon \oplus \cdots \oplus \epsilon}_n
$$
\vskip -1.2cm
\hfill {$(\Gamma)$}
\vskip0.6cm
\noindent
over the compact space $\RR P^1 \simeq S^1$.

We indicate several proofs that such an isomorphism is impossible, starting from the ``mechanical'' point of view. The method of characteristic has a dual nature relative to our previous ``mechanism'', which was based on the Steenrod cohomology operations. To carry this approach through, it is enough to fall  back on definitions and axioms.

For any real vector bundle $\xi$ over $B$. one has an $i$-th Stiefel-Whitney class  $w_i (\xi) \in H^i (B, \ZZ_2)$. It is also axiomatic that the canonical line bundle $\gamma$ over $\RR P^1$ is a non-zero element in  $H^1(\SS^1, \ZZ_2)$, the only one. For our simple purposes, we work with the {\it total} Stiefel-Whitney class which is defined as
$$w(\xi) = 1 + w_1 (\xi) + w_2(\xi) + \cdots + w_k (\xi)\, + \,0\, \, \cdots \,,$$
where $k$ is taken as the dimension of the CW-complex  or manifold $B$. In this formulation we get to deal with larger bundles through the Whitney Product Formula $w(\xi \oplus \eta) = w(\xi) \cup w(\eta)$. In general where characteristic classes are concerned, ``product'' means ``cohomology cup product''. Also $w(\epsilon) = 1$  by \cite{Milnor \&  Stasheff}, so our bundle isomorphism $(\Gamma)$ turns into
$$\left[1 + a\right]^n = \roman 1 \in H^* \left(\SS^1, \ZZ_2\right)\,,$$
the latter seen more generally as the graded ring  $\prod_{j=1}^{\infty}\,H^i (B, \ZZ_2)$. Here $w_1 (\gamma) = a$. But since $a \cup a = 0$ in $H^2(\SS^1, \ZZ_2)$, a calculation gives $w(n \,\gamma) = 1 + a \ne 1$, when $n$ is odd.

Thus the two $n$-plane bundles, proposed to be equal, are nothing of the sort. Hence two such $n\times n$ matrices (``incomparable'' by scale) do not exist and any such real matrix has a real eigenvalue and eigenvector. A similar ``mechanical'' proof will be applied to the  less obvious case, namely where $q = 2n$. We stick with the $q = \text{odd}$ case for a moment, to illustrate the geometric or ``handmade'' point of view.

Let $\text{Vect}_n (B)$ denote the set of $B$-ismorphism classes of $n$-plane bundles. By the use of clutching functions it is not hard to characterize this set completely in case $B = \SS^1$, \cite{Dugger, p. 68}. In fact the only (two) elements are $\epsilon^n$ and $\gamma \oplus \epsilon^{n-1}$. This classification immediately shows that $(\Gamma)$ doesn't work, and we recover our real eigenvalue and eigenvector, for a given $A_{n\times n}$.
Finally, we try some algebra along even more elementary lines. We agree that mathematical civilization has stumbled upon Whitney sums and the canonical line bundle, but the discovery of cellular cochains and their boundary-equivalence classes awaits a remote future era.

We note at once that $\gamma \oplus \gamma \simeq \epsilon \oplus \epsilon$, that is to say  that $2\gamma$ is ``trivial''. To show the one has to  come up with two linearly independent ``sections'' on $\SS^1$. Using a standard Cartesian embedding $\SS^1 \to \RR^2$, we propose  $\left[\matrix
x & y \\
-y & x
\endmatrix\right]$ as the sections, one to each row. These do not vanish on $\SS^1$ in view of $x^2+y^2=1$, and are independent by virtue of the inner product. Most importantly, the sections are antipode preserving (equivariant relative to $\pm 1$ so they lead to sections $\sigma_j = \RR P^1\to \RR^2,\, j = 1,2$. We have $\sigma_j (-s) = -\sigma_j(s)$, $s \in \RR P^1$. On the other hand, $\gamma \oplus \epsilon^{n-1}$ cannot be trivial (also for $n$ even). For since $\text{dim}\,\SS^1  = 1$, Proposition 8 shows that $\gamma \oplus \epsilon$ would be trivial too. But imagine sections
$$\matrix
x_1(s) & y_1(s) \\
x_2(s) & y_2(s)
\endmatrix = \left[\Lambda(s)\right]\,\,,$$
satisfying $x_1(-s) = -x_1(s)$, (component of Hopf line) $y_1(-s) = +y_1(s)$, (component of the trivial line), and similarly $x_2(-s) = -x_2(s)$, $y_2(-s) = y_2(s)$. Then  if $\text{det}\,\Lambda (s_0) = \lambda \ne 0$, we would have $\text{det}\,\Lambda (-s_0) = -\lambda$. Continuity of the sections, together with connectedness of the circle and the Intermediate Value Theorem, shows that $\Lambda (t)$ has deficient rank at some $t \in \RR P^1$. Hence the sections were not truly independent. We come to the same conclusion as before, that $n \gamma$ is {\it not} trivial when $q = n = 2m+1$. This last argument is what we prefer, to prove that $A_{n\times n}$ has a real eigenvalue-eigenvector pair.

We now move to the less trivial case where $q = 2n$, the bases space of bundles $\gamma$ and $\epsilon$ is $\RR P^2$ and our three independent, nonsingular and spanning matrices $A_1, A_2, A_3$ are of size $  q\times q$. The alleged isomorphism that  arises is
$$q\,\gamma \simeq \epsilon^q \ \,.$$
\vskip -0.8cm
\hfill {$(\Delta)$}
\vskip0.4cm
As in the previous one-dimensional instances, we may apply the total Stiefel-Whitney operation to both sides of the isomorphism. Can $w(q\,\cdot \gamma)$ equal $1$? Recall that if $w(\gamma) = 1+a$, we have a defining ring relation $a^3 = 0$ (cup product) in $H^* (\RR P^2, \ZZ_2)$.
Jotting down the ($\text{mod}\,2$) Pascal triangle to compute $w(k\,\gamma)$ leads to $1 + a^2 \ne 1$ when $k$ is twice an odd integer. Thus $(\Delta)$ is refuted and every  complex square matrix of dimension $n = 2m+1$ has an eigenvector.

We felt that the use of characteristic classes and cohomology was a step away from the basic geometry that lies behind

\bigskip
\noindent
{\bf Theorem.} Any three real square matrices $A_1, A_2, A_3$ of size $q = 2n = 2(2m+1)$ can be combined non-trivially to give a singular matrix.

\bigskip
Consider the approach of classifying $\text{Vect}_1 (\RR P^2), \text{Vect}_2(\RR P^2), \dotsc$, similar to \linebreak what we did with $\text{Vect}_k(\SS^1¡)$. The study of projective space rather than the sphere is more complicated. Part of the answer lies in the  bijection $\text{Vect}_n(\SS^2) \simeq \pi_1 (SO_n)$. The latter fundamental group is trivial when $n = 1$, and equals $\ZZ_2$ for $n > 2$. More generally, a ring structure induced by {\it tensor product} of vector bundles can be imposed onto the Grothendieck (additive) group of ``stable isomorphism classes'' of vector bundles. This ring structure on $\tilde{K}_{\RR}(\RR P^{r-1})$ is probably more than what is needed for issues concerning {\it constant rank} matrices. To compute with this ring brings the ``Atiyah-Hirzebruch spectral sequence'' into play. Such powerful techniques leave essential questions open, for example: is $\gamma \oplus \gamma$ stably trivial over $\RR P^2$? We see below that it is not, by showing that $\gamma \oplus \gamma \oplus \epsilon$ is non-trivial (using Propositions 7 \& 8). Another  argument is sketched out in \cite{Dugger, p. 82},  observing that the self-intersection of an $\RR P^1$ embedded in $\RR P^2$ contains a generic point.

Finally we arrive at a straight homotopical  analysis of equation $(\Delta)$, which constitutes whatever  innovation this report may offer. If we see that $(\Delta)$ cannot be valid, we obtain equivalent to the Theorem above: any complex linear transformation $T_{n\times n}$, $n$ odd, has a (complex) eigenvalue-eigenvector pair.
Similar to how we handled the simple case of equation $(\Gamma)$, we break off ``trivial'' chunks of the left-hand side. In fact,

\bigskip
\noindent
{\bf Lemma.} Over $B = \RR P^2$ we have $\gamma \oplus \gamma \oplus \gamma \oplus \gamma  = \epsilon \oplus \epsilon \oplus \epsilon \oplus \epsilon$. That is, $4\gamma \simeq \epsilon^4$.

\bigskip
\noindent
{\it Proof.} Take $\SS^2$ embedded as usual in $3$-space. Then we look for four independent sections. But these may be racked up as rows in a matrix
$$\matrix
x & y & z & 0 \\
-y & x & 0 & z \\
-z & 0 & x & -y \\
0 & -z & y & x
\endmatrix \ \ \,\,\,\,.
$$
\indent
A row vector as a function on $\SS^2$ never vanishes in view of $x^2+y^2+z^2=1$. All distinct pairs of rows are orthogonal, hence the whole set of rows is independent. Or we could calculated the determinant of this $4 \times 4$ array, yielding the non-zero constant $1$ on $\SS^2$. The significant observation is that each given section $\sigma$ satisfies $\sigma (-s) = -\sigma (s)$ so it represents a ``canonical'' section on every component.\qed

Since $q = 2n = 2(2m+1)$ we end up with
$$\gamma \oplus \gamma \oplus \epsilon^{4m} \simeq \epsilon^{4m+2}\,\,.$$
As $\RR P^2$ has dimension 2, Propositions 7 \& 8 yield
$${\Cal X} := \gamma \oplus \gamma \oplus \epsilon \simeq \epsilon^3\,\,.$$
The remainder of the article shows that this bundle isomorphism is absurd. But Dugger's  argument could also be applied, and in any case should be generalized. This would prove our  Theorem by elementary intersection theory.

As it is we use two geometric facts, both  applied to the 2-sphere. The first of  them is the well-known Borsuk-Ulam Theorem. We emphasize that  our understanding of complex linear mappings of  any odd order, however large, boils down to  topological observations regarding $\SS^2$ only. The other fact is that an antipode-preserving  mapping $g: \SS^2 \to \SS^2$ is never homotopic to  its antipode $-g$. This is because for $\SS^2 \subset \RR^3$, negation switches orientation, hence changes the sign of the mapping degree, which was non-zero to start with in view of B-U Theorem.

Various equivalent statements of the B-U Theorem, with proofs, can be found in \cite{Matou$\check{\text{s}}$ek}, in standard topology texts, and in monographs available on-line. We emphasize one formulation, ``there is no mapping $\psi: \SS^2 \to \SS^1$ that {\it preserves antipodes}'', namely $\psi(-s) = -\psi (s)$, where $\SS^1$ could be the equator of $\SS^2$, the locus of {$z = 0$}. In other words, there is no mapping $\psi_*: \RR P^2 \to \RR P^1$, period. ``Degree Theory'' is one  of our covert themes. We review the fact that an antipode-preserving $\phi: \SS^1 \to \SS^1$ has {\it odd degree} and therefore cannot be extended to a disk. 

One may decompose $\SS^1$ into an upper arc $\tau$ and a lower arc $\sigma$. These arcs are seen as mappings from  $[0,1]$, covering $\SS^1$ from $(1,0)$ to $(-1,0)$ to $(1,0)$ in a counterclockwise manner. In fact we can take  $\tau (t) = \left(\cos \pi t, \sin \pi t\right)$ etc. The book of \cite{Sieradski} points out that $\varphi \tau (1) = \varphi \tau (0) + (2m+1) \pi$ for some integer $m$. That is, the total length of $\varphi(\tau)$, with reference to the ``fundamental groupoid'', equals $L$, an odd integer times $180^{\circ}$. But $\sigma$, the lower hemisphere arc on $\SS^1$ that matches $\tau$ by orientation, must also yield total $\text{length}[\varphi(\sigma)] = L$. Thus $\text{length}\,\varphi(\SS^1) = \text{length}\,\varphi(\tau) + \text{length}\,\varphi(\sigma) = 2L = \text{odd multiple of $2\pi$}$, hence non-zero. The theory of the fundamental groupoid, which uses some linear approximation, shows that a contractible $\varphi: \SS^1 \to \SS^1$ (extendible to the bounded disk) has  ``winding number $0$'' in this sense.

Note once more that we invoke only the  B-U Theorem on $\SS^2$, another equivalent wording being ``any mapping $\SS^2 \to \RR^2$ takes some pair of antipodal points to the {\it same} point.'' This is the first substantial case of the B-U Theorem; higher-dimensional versions are more complicated to prove.

While we are at it, let us mention  a related analytic demonstration for $\SS^2$ due to A.  Carbery. He also builds up a proof of B-U for  dimension $n$, namely that ``a mapping $\kappa: \SS^n \to \RR^n\,$, equivariant with respect to negation, has a zero $s_0 \in \SS^n$''; that is, $\kappa (s_0) = \vec{0} \in \RR^n$.

Suppose we have $g = (g_1, g_2) = \DD^2 \to \SS^2$, equivariant on the boundary. That is, $g(-t) = -g(t)$ for $t \in \SS^1$. Now the Jacobian matrix $J_g$ has rank at most $1$. We apply Stokes' theorem to yield
$$0 = \int_{\DD^2} dg = \int_{\DD^2} \text{det}\,J_g\,dx \wedge dy = \int_{\DD^2} dg_1 \wedge dg_2 = \int_{\SS'} g_1\,dg_2\,\,,$$
so that if we show that
$$\int_{\SS^1} g_1\,dg_2 \ne 0\,,$$
we would have a contradiction. Even better, the same formulation also yields $0 = -\int g_2\,dg_1$. It would now be enough to show that
$$\int_0^1 \left[g_1(t) g_2'(t) - g_2 (t) g_1' (t)\right] dt \ne 0\,.$$
But this integral amounts to $\int_{\SS^1} d \theta$, where $\theta$ gives counter-clockwise oriented arc length, again the antipode preservation implies that this arc length equals twice
$\int_0^{1/2} ds(t) = $ length from $\varphi(0)$ to $\varphi (\frac12) \ne \varphi (0)$, so is non-zero.

Finally, it is hard to improve upon a modern ``classical'' treatment of the Borsuk-Ulam Theorem as found in \cite{Massey}, where the author works with the fundamental group and its  action on covering spaces such as $\SS^1 \to \RR P^1$ and $\SS^2 \to \RR P^2$.

Returning to the main issue, if ${\Cal X} = \gamma \oplus \gamma \oplus \epsilon$ were a trivial bundle, there would be three independent sections $\tau_j : \RR P^2 \to \RR^3$. The first two components of these sections could be represented by antipode-preserving maps with $\tau_j^i(-s) = -\tau_j^i(s)$ for $i = 1,2$. Each {\it third} component represents mappings into a trivial line, so that
$$\tau_j^3(-s) = \tau_j^3(s), \quad s \in \SS^2\,, \;j = 1, 2, 3\ \,.$$
\indent
In particular, the three-by-three matrix function $\Omega_{ij}(s) = \left[\tau_i^j(s)\right]$ must be nonsingular for each $s \in \SS^2$. We demonstrate to the contrary that a rank reduction of  this matrix must occur somewhere. We start with a new geometric interpretation of the first two columns. The first and second columns are really mappings  $\rho_k:  \SS^2 \to \RR^3$, $k = 1,2,\,$ which looked at coordinate-wise, turn out to be antipode preserving. In case $\rho_1$ ever achieves the  origin, $\rho_1(s_0) = (0,0,0)$, $\text{det}\,\Omega(s_0)$ goes to $0$ immediately. Thus, the $\{\tau\}$ are not the trivializing sections that we thought they were. We may assume that $\rho_1, \rho_2$ avoid the value $\vec{0} \in \RR^3$. Normalizing to $\dfrac{\rho(s)}{\|\rho(s)\|}$ but keeping the notation, we obtain two antipode-preserving
$$\rho_j: \SS^2 \to \SS^2\ .$$
By the 2-dimensional Borsuk-Ulam Theorem, both of them are {\it onto} maps, homotopically essential, and have {\it odd} Brouwer degree.

Furthermore, at a given $s \in \SS^2$, $\rho_1$ never takes the same value as does $\rho_2$, lest the matrix $\Omega(s)$ lose rank. Elementary homotopy theory, see \cite{Dugundji} shows that $\rho_1 \sim \text{anti}\,\rho_2$ (they are freely homotopic maps $\SS^2 \to \SS^2$). For the same reason, $\rho_1(s)$ can never equal $-\rho_2(s)$ either, hence $\rho_1 \sim \rho_2$. We conclude that $\rho_2$ is an onto mapping with odd (non-zero) degree and $\rho_2 \sim -\rho_2$. But for a self-map $g$ of a manifold in {\it even} dimensions, one has  $\text{deg}(-g) = -\text{deg}(g)$, so a contradiction follows from the homotopical invariance of the Brouwer degree.

In carrying out the definition of Brouwer degree, homotopical invariance and other essential features, there are several popular approaches. Since we have chosen to stand away from the  selection of special points in a space, which uses Sard's theorem. Thus we leave for the reader to carry out our Theorem by means of such techniques from oriented intersection theory/regular values, \cite{Milnor TFDV, Guillemin \& Pollack}.

A straightforward alternative involves simplicial approximation, where the degree as an integer emerges immediately, and the homotopy invariance is natural. Returning to the question of whether $g \sim -g$ on $\SS^2$ could occur, we have available the Hopf-Freudenthal theory, which characterizes such mappings up to homotopy (as the suspension of some $f: \SS^1 \to \SS^1$). The problem reduces to showing that $\text{id}_{\,\SS^n}$ cannot be homotopic (on an even sphere) to $-\,\text{id}_{\,\SS^n}$, the antipoda1 map. This result is {\it stronger} than the actual Borsuk-Ulam result. Indeed, from the standpoint of differential topology, the latter can be proved by modulo-two methods, unlike the former.

Whether we wish to use Hopf-Freudenthal or not, the validity of integer degree shows that  our matrix $\Omega(s)$ must somewhere be singular on $\SS^2$, giving that ${\Cal X}$ is a non-trivial bundle over $\RR P^2$, that bundle ``equation'' $(\Gamma)$ does not hold, and that any $A: \CC^n \to \CC^n$ linear, has an eigenvector when $n$ is odd.

The analytic definition of ``degree'' exists in various versions, for instance the Heinz-Mawhin degree which enjoys continued research. Homotopical invariance comes easily, \cite{Dinca-Mawhin}. This definition has solid historical roots in its relation to Kronecker's characteristic or ``curvature integra'', which goes back at least to  Gauss. Proving that the degree for continuous maps is an integer was the arduous aspect, but new proofs have emerged. Besides, we only need to know that antipode-preserving $g$ satisfies $\text{deg}(g) \ne 0$, and that  $\text{deg}(-g) = -\text{deg}(g)$, which come out of the defining formulas \cite{Heinz}.

By either of these methods we have proven that $\gamma \oplus \gamma$ is not a stably trivial plane bundle over $\RR P^2$. It was previously noted by D. Dugger that  the spectral sequence associated to cellular decomposition has been used to prove this. He  suggested a more natural geometric approach to  this type of question. We hope that suitable insight and extension of the method presented here, would lead to alternative K-theory calculations and even new results on the frontier of  topology and linear algebra.

\end{document}